\documentclass[a4paper,twocolumn,colorlinks]{preprint}

\usepackage[english]{babel}
\usepackage{graphicx}
\graphicspath{{figures}}

\usepackage{amsmath}
\usepackage{amsthm}
\usepackage{amssymb}
\usepackage{mathtools}
\usepackage{enumitem}
\usepackage{siunitx}

\usepackage{tikz}
\usepackage{pgfplots}
\pgfplotsset{compat=newest}
\usetikzlibrary{3d}
\usetikzlibrary{decorations.shapes}
\usetikzlibrary{shapes.arrows}
\usetikzlibrary{calc}
\usetikzlibrary{arrows}

\def\n{2.2}
\def\p{.5}
\def\m{1}
\def\gap{0.025}
\def\sat{30}
\def\opac{1}
\newcommand\snapshot[2]{
  \begin{scope}[canvas is yz plane at x=#1]
    \draw[fill=red!\sat,opacity=\opac,cm={1,0,0,1,(-\gap cm,\gap cm)}] (0,\p) rectangle (\n,\n+\p) node[pos=0.5,transform shape,red!\sat!black] {\(A_{p_{#2}}\)};
    \draw[fill=blue!\sat,opacity=\opac,cm={1,0,0,1,(\gap cm,\gap cm)}] (\n,\p) rectangle (\n+\m,\n+\p) node[pos=0.5,transform shape,blue!\sat!black] {\(B_{p_{#2}}\)};
    \draw[fill=green!\sat,opacity=\opac,cm={1,0,0,1,(-\gap cm,-\gap cm)}] (0,0) rectangle (\n,\p) node[pos=0.5,transform shape,green!\sat!black] {\(C_{p_{#2}}\)};
    \draw[fill=orange!\sat,opacity=\opac,cm={1,0,0,1,(\gap cm,-\gap cm)}] (\n,0) rectangle (\n+\m,\p) node[pos=0.5,transform shape,orange!\sat!black] {\(D_{p_{#2}}\)};
  \end{scope}
}

\tikzset{decorate sep/.style 2 args=
  {decorate,decoration={shape backgrounds,shape=circle,shape size=#1,shape sep=#2}}}
\usepackage[linesnumbered,ruled,algo2e,norelsize]{algorithm2e}
\SetKwInOut{Input}{Input}
\SetKwInOut{Output}{Output}

\usepackage{bbm}

\usepackage{mathrsfs}

\AddToHook{begindocument/end}{
\theoremstyle{plain}
\newtheorem{theorem}{Theorem}[section]
\newtheorem{lemma}[theorem]{Lemma}

\theoremstyle{definition}
\newtheorem{definition}{Definition}[section]
\newtheorem{example}{Example}[section]
\newtheorem{remark}[theorem]{Remark}}

\usepackage{etoolbox}
\AtBeginEnvironment{bmatrix}{\setlength{\arraycolsep}{.27em}}

\begin{document}


\title{Snapshot-driven Rational Interpolation \texorpdfstring{\\}{}of Parametric Systems}

\author[$\ast$]{Art J. R. Pelling}
\affil[$\ast$]{Department of Engineering Acoustics, Technische Universität Berlin, Einsteinufer 25, 10587 Berlin, Germany.\authorcr
  \email{a.pelling@tu-berlin.de}, \orcid{0000-0003-3228-6069}}

\author[$\dagger$]{Karim Cherifi}
\affil[$\dagger$]{Department of Mathematics, Technische Universität Berlin, Straße des 17. Juni 136, 10623 Berlin, Germany.\authorcr
  \email{cherifi@math.tu-berlin.de}, \orcid{0000-0003-1294-9291}}

\author[$\ddagger$]{Ion Victor Gosea}
\affil[$\ddagger$]{Max Planck Institute for Dynamics of Complex Technical Systems, Sandtorstr. 1, 39106 Magdeburg,
Germany.\authorcr
  \email{gosea@mpi-magdeburg.mpg.de}, \orcid{0000-0003-3580-4116}}

\author[$\S$]{Ennes Sarradj}
\affil[$\S$]{Department of Engineering Acoustics, Technische Universität Berlin, Einsteinufer 25, 10587 Berlin, Germany.\authorcr
  \email{ennes.sarradj@tu-berlin.de}, \orcid{0000-0002-0274-8456}}

\shorttitle{Snapshot-driven rational interpolation of parametric systems}
\shortauthor{A. J. R. Pelling, K. Cherifi, I. V. Gosea, E. Sarradj}

\abstract{
Parametric data-driven modeling is relevant for many applications in which the model depends on parameters that can potentially vary in both space and time. In this paper, we present a method to obtain a global parametric model based on snapshots of the parameter space. The parameter snapshots are interpolated using the classical univariate Loewner framework and the global bivariate transfer function is extracted using a linear fractional transformation (LFT). Rank bounds for the minimal order of the global realization are also derived. The results are supported by various numerical examples.}


\novelty{A new interpolation method based on the Loewner framework and a linear fractional transformation for parametric system realizations is presented. It can be used in existing matrix interpolation frameworks for parametric linear time-invariant systems.}


\maketitle

\section{Introduction}

Many dynamical systems encountered in science and engineering problems \cite{bock2013,lu2017,vojkovic2022} exhibit parameter-varying dynamics. Instances of dynamics-influencing parameters can be encountered when considering material properties in mechanical systems, wing shape, and weight distribution in aircraft design, varying source and receiver locations in acoustical transmission, and variable resistors in circuit design.

In the linear case, these parametric systems can be expressed in the form of a parametric linear time-invariant (pLTI) systems of the form:
\begin{align}\label{eq:DAE}
\begin{split}
		\tfrac{d}{dt} x(t,p) &= A(p)x(t,p) + B(p)u(t),\\ 
	   y(t,p) &=C(p)x(t,p) + D(p)u(t),
	\end{split}
\end{align}
with the so-called \emph{transfer function} that is obtained by a Laplace transform
\begin{equation}
    \label{eq:ptf}
    H(s,p)=C(p)(sI-A(p))^{-1}B(p)+D(p),
\end{equation}
where $u: \mathbb{R} \to \mathbb{R}^{n_{\mathrm{i}}}$, $x: \mathbb{R}\times P \to \mathbb{R}^n$, $y : \mathbb{R}\times P \to \mathbb{R}^{n_{\mathrm{o}}}$ are the \emph{input}, \emph{state}, and \emph{output} of the system, respectively, and $A(p)\in\mathbb{R}^{n\times n}, B(p)\in\mathbb{R}^{n\times n_{\mathrm{i}}}$, $C(p)\in\mathbb{R}^{n_{\mathrm{o}}\times n}$, $D(p)\in\mathbb{R}^{n_{\mathrm{o}}\times n_{\mathrm{i}}}$, and $p\in P\subset\mathbb{R}$ is a scalar parameter.

In some applications \cite{kujawski2024,lu2017,vojkovic2022}, a general system expressed as \eqref{eq:DAE} is not readily available in the sense that the model might only be available for a few parameter values $p_i$ because the system is either modeled from data where only a discrete sampling of the parameter range is possible or the order of the original system is simply too large and it is necessary to obtain a reduced order global model. This problem led to an increased interest in developing data-driven modeling and model reduction techniques for such parametric systems. An in-depth comprehensive survey is available in \cite{benner2015} which contains a plethora of interpolation-based reduction methods for pLTI systems.


Our work focuses on the data-driven interpolation of the system in \eqref{eq:DAE}, i.e. an interpolation of local non-parametric models according to \cite[Sec. 4.2]{benner2015} where each non-parametric model corresponds to a fixed parameter value \(p_i\). This field can be divided into three different categories that differ with regard to the input data and model premises. The classes of methods that interpolate local bases will be disregarded here due to unfavorable computational properties described in \cite[Sec. 4.2.1]{benner2015}. Instead, we will introduce the two remaining classes, namely \emph{matrix interpolation} and \emph{transfer function interpolation} below. So-called matrix interpolation denotes the class of methods that interpolate local (reduced) model matrices \cite[Sec. 4.2.2]{benner2015}. Our proposed method belongs to this class. Additionally, we will also briefly outline the interpolation of local (reduced) transfer functions \cite[Sec. 4.2.3]{benner2015}, since methods of this class are always admissible under the prerequisites of our proposed method. We now proceed with a short description of both method classes:

\begin{description}[leftmargin=0cm]
    \item[Matrix Interpolation]
    Given snapshots of realizations
    \begin{equation}
    \label{eq:snapshots}
        \mathbf{\Sigma}_i=\bigl(A(p_i),B(p_i),C(p_i),D(p_i)\bigr),
    \end{equation}
    of the system matrices at $n_{p}\in\mathbb{N}$ fixed parameter samples \(p_i\) and denote
    \begin{equation*}
        \Pi:=\left\{p_i~|~i=1,\,\dots,\,n_p\right\}\subset P
    \end{equation*}
    the set of all parameter samples. The aim is to construct a parametric transfer function $\hat{H}(s,p)$ that interpolates the data, i.e. \(\hat{H}(s,p_i)=H(s,p_i)\), for any given value of $p_i\in\Pi$ and $s\in\mathbb{C}$.
    The interpolation of the system matrices is straightforward and standard matrix interpolation techniques can be employed. Many works \cite{panzer2010,geuss2013,geuss2014,vojkovic2022,yue2019} use direct interpolation, such as
    \begin{equation}
        \label{eq:matrix-interpolation}
        \hat{\mathcal{M}}(p)=\sum_{i=1}^{n_{p}}\omega_i(p)\mathcal{M}(p_i),\quad\mathcal{M}\in\{A,B,C,D\}
    \end{equation}
    where \(\mathcal{M}\) is any of the matrix functions in \eqref{eq:DAE} and \(\omega_i:P\rightarrow\mathbb{R}\) is a weighting function that fulfills \(\omega_i(p_j)=\delta_{ij}\). Optionally, matrix manifold interpolation \cite{amsallem2011,geuss2013} can be used, where the interpolation is in the tangent space of a matrix manifold instead of the matrices directly. Using the \emph{Riemannian logarithm} and \emph{Riemannian exponential}, one can locally map the tangent space and manifold onto each other. For an introduction to parametric model order reduction with matrix manifold interpolation, see \cite{zimmermann2021}.
    
    Unfortunately, the interpolation of the snapshots in \eqref{eq:snapshots} at the parameter samples \(p_i\) alone does not imply that the interpolated model matches the true system at intermediate parameter values \(p\in P\setminus\Pi\). The approximation quality of the interpolated models can be assessed in various ways such as computing the differences in the entries of the realization matrices, as done in \cite{amsallem2011}. However, a global quality measure in the form of multivariate transfer function norm, e.g. the \(\mathcal{H}_2\otimes\mathcal{L}_2\) norm \cite{hund2022,baur2011}, is preferable. A common observation is that the modeling error at intermediate parameter values only becomes small if all of the snapshots \eqref{eq:snapshots} are expressed in common state coordinates \cite{amsallem2011}. If the snapshots stem from explicit equations or computer simulations this requirement is usually met. However, in these cases, the snapshot matrices may be high-dimensional which leads to an increased cost of evaluating \eqref{eq:matrix-interpolation} and consequently \(\hat{H}(s,p)\). The computational performance of the parametric model can be greatly improved by employing locally reduced order system matrices of order $r$ as input snapshots, i.e.
    \begin{equation}
        \label{eq:rommats}
    \begin{aligned}
\mathbf{\Sigma}_{i,r}&=\bigl(A_{r,i},B_{r,i},C_{r,i},D_{r,i}\bigr)\\&=\left(T_{1,i}^TA(p_i)T_{2,i},T_{1,i}^TB(p_i),C(p_i)T_{2,i},D(p_i)\right),
    \end{aligned}
    \end{equation}
    with projection matrices \(T_{1,i},T_{2,i}\in\mathbb{R}^{n\times r}\) such that the condition \(T_{1,i}^TT_{2,i}=I_r\) holds. The projection matrices originate from a projection-based model order reduction algorithm of choice. In general, these reduced order non-parametric models will no longer share a common state space, when being reduced individually. Initially suggested in \cite{panzer2010}, a reprojection of the reduced order system matrices onto a common generalized reduced order coordinate system is necessary here in order to achieve meaningful models for intermediate values (see also \cite{amsallem2011,geuss2013}).
    
    A general framework for matrix interpolation is offered in \cite{geuss2013}. The framework unfolds the different algorithmic variations of the previous works and includes manifold interpolation and model reduction with different reprojection techniques as optional steps. We consider our work to be in line with this framework and contribute an alternative interpolation method using rational interpolation with the Loewner framework \cite{mayo2007}. For the sake of brevity, we do not explicitly consider manifold interpolation nor local model order reduction as they can be easily incorporated by following the pertaining steps in \cite{geuss2013}.
\begin{remark}
\label{rem:descriptor}
The general framework in \cite{geuss2013} considers systems in \emph{descriptor form} \((E(p),A(p),B(p),C(p))\) instead of the \emph{`standard realization`} in \eqref{eq:DAE}. Nevertheless, this is by no means a restriction, as these two model structures can be transformed into one another under certain conditions, see \cite{kunkel2006} for details. 
\end{remark}
    \item[Transfer Function Interpolation]
    Naturally, one can also interpolate the bivariate rational transfer function \(H(s,p)\) directly. Instead of snapshots of the underlying system matrices, samples of the transfer function \(H(s_i,p_i)\) at given frequency-parameter tuples \((s_i,p_i) \in\mathbb{C}\times\mathbb{R}\) can be used for interpolation via
    \begin{equation*}
        \hat{H}(s,p)=\sum_{i=1}^n\omega_i(s,p)H(s_i,p_i),
    \end{equation*}
    with weighting functions \(\omega_i:\mathbb{C}\times P\rightarrow\mathbb{R}\). The specific choice of weighting function might impose restrictions on the sampling strategy, e.g. \((s_i,p_i)\) lying on a grid. An overview of polynomial bivariate interpolation can be found in \cite{olver2006}. If one has access to an underlying realization, the concept of model reduction can also be applied. In \cite{baur2009,baur2011a}, samples of reduced-order transfer functions are interpolated with polynomial and rational interpolation techniques, respectively. 

Another approach is multivariate rational approximation, for which the Loewner framework is a prominent representative in the field of system and control theory and can be directly used for transfer function interpolation and model reduction purposes.
Although originally developed for the univariate case \cite{mayo2007}, the method was extended to the bivariate \cite{lefteriu2011,antoulas2012} and multivariate case \cite{ionita2014}. The capabilities of this approach concerning model order reduction of parametric dynamical systems are investigated in \cite{lefteriu2011,ionita2014}. At the core of this approach stands the multivariate Loewner matrix whose rank encodes information about the minimal order of the approximating multivariate rational function. This is precisely what sets the approach apart from other, more conventional tools for polynomial and rational approximation.
Additionally, \cite{antoulas2012} provides an effective way to construct realizations of linear parametric models directly from the data (samples of the multivariate transfer functions). More recently, \cite{vojkovic2023} proposed an efficient way to deal with the case of multiple inputs and multiple outputs.

The curse of dimensionality (C-o-D) may limit applications of the Loewner framework with an increased number of parameters. However, as shown in the recent review paper \cite{antoulas2024}, the C-o-D may be tamed by alleviating the explicit computation of large-scale multivariate Loewner matrices.
  Furthermore, another potential bottleneck of the Loewner framework is that the (bivariate) Loewner matrix depends on a particular choice of partitioning for samples \((s_i,p_i)\). To date, it is not understood how to choose optimal partitionings (especially in the parametric case), since this strongly depends on the application under investigation. One \emph{suitable} partitioning scheme was recently proposed in the pAAA (parametric Adaptive Antoulas-Anderson) algorithm  \cite{rodriguez2023}. There, similarly to the original AAA algorithm \cite{nakatsukasa2018}, a greedy scheme of selection is proposed, and the dimension of the rational approximant is increased at each step until a desired approximation quality on the data is reached.

The method proposed in this paper does not involve multivariate Loewner matrices, and hence it is not directly related to the parametric Loewner framework, as far as the authors are aware.
Consequently, in the current study, we harness the simplicity, elegance, and effectiveness of the standard, univariate Loewner framework and show how can this be used for specific parametric problems.
\end{description}
\begin{remark}
    A key difference of this type of interpolation problem as opposed to the matrix interpolation problem mentioned earlier is that the sampling data is assumed to be known only at a discrete set of frequencies \(s_i\in\mathbb{C}\), whereas \eqref{eq:ptf} offers a way to directly evaluate \(H(s,p_i)\) over the entire frequency range at \(p_i\) from the snapshots. Hence, by choosing a set of frequencies \(s_i\in\mathbb{C}\), the necessary sampling data for transfer function interpolation algorithms can always be computed from the system snapshots in \eqref{eq:snapshots}.
\end{remark}
\begin{sloppypar}
There exist numerous other methods, e.g. 
\cite{aumann2021,benner2014,benner2024,feng2019,feng2024,gosea2021,goyal2024,grimm2018,hund2022,mlinaric2023}, that are out of the scope of this paper and will not be discussed further.
\end{sloppypar}
The objective of this work is to present a data-driven realization technique to obtain a bivariate parametric transfer function from system snapshots using the Loewner framework. After introducing some preliminaries and the Loewner framework in \cref{sec:Preliminaries}, we present the derivation and procedure of our method in \cref{sec:method}. In \cref{sec:rank}, we present rank bounds on the obtained Loewner matrices and by consequence the order of the minimal realization of the global parametric realization. The theoretical results are supported by numerical examples in \cref{sec:numerics} where we use our algorithm to find the parametric transfer function of different parametric model classes.

\section{Preliminaries}
\label{sec:Preliminaries}
The method presented in this paper relies on univariate Loewner interpolation of the parametric model and uses a \emph{linear fractional transformation} (LFT) for obtaining a general transfer function. In this section, an introduction to these two concepts is presented. 

\subsection{Linear fractional transformation}
\label{sec:lft}
For notational convenience, we will represent transfer functions by means of an LFT as suggested in \cite{doyle1991}.
\begin{definition}[\cite{doyle1991}]
\label{def:LFT}
Let \(S\) be a complex matrix partitioned as
\begin{equation*}
    S=
    \begin{bmatrix}
    S_{11}&S_{12}\\S_{21}&S_{22}
    \end{bmatrix}\in\mathbb{C}^{(m_1+m_2)\times(n_1+n_2)}.
\end{equation*}
Let \(\Delta\in D\subset\mathbb{C}^{n_1\times m_1}\) with \((I_{m_1}-S_{11}\Delta)\)  non-singular. Then, we define the map \(\mathcal{F}_{u}(S,\,\cdot\,):D\rightarrow\mathbb{C}g^{m_2\times n_2}\) as
\begin{equation*}
    \mathcal{F}_{u}(S,\Delta)=S_{21}\Delta(I_{m_1}-S_{11}\Delta)^{-1}S_{12}+S_{22}
\end{equation*}
as the \emph{(upper) linear fractional transformation} (LFT). If \(m_1=n_1\) and \(\Delta\) is non-singular, we can also write it as
\begin{equation*}
    \mathcal{F}_{u}(S,\Delta)=S_{21}(\Delta^{-1}-S_{11})^{-1}S_{12}+S_{22}.
\end{equation*}
\end{definition}
The transfer function of an LTI system \(\mathbf{\Sigma}=(A,B,C,D)\) can be expressed by the following LFT for \(s\neq0\)
{\delimitershortfall=0pt
\begin{equation}
\label{eq:LFT-tf}
G(s)=C(sI_n-A)^{-1}B+D=\mathcal{F}_{u}\left(\begin{bsmallmatrix}A&B\\C&D\end{bsmallmatrix},s^{-1}I_n\right).
\end{equation}
}
The representation of the transfer function of an LTI system as an LFT of a larger matrix is what enables the interpolation with the univariate Loewner framework and ultimately leads to the main result of this paper in \cref{sec:method}.
\begin{remark}
\label{rem:schur}
It should be noted that \eqref{eq:LFT-tf} is only well-defined for \(s\neq0\). However, for \(s=0\), the transfer function can be written as the Schur complement
{\delimitershortfall=0pt
\begin{equation}
\label{eq:schur}
    G(0)=D-CA^{-1}B=\begin{bsmallmatrix}
        A&B\\C&D
    \end{bsmallmatrix}/A.
\end{equation}
}
The distinction between these two cases is a technicality that arises from \cref{def:LFT} and will not become relevant to the contents of this work. Hence, by abuse of notation, we will implicitly refer to \eqref{eq:schur} by the LFT in \eqref{eq:LFT-tf} whenever \(s\) is zero.
\end{remark}

\subsection{Loewner matrix interpolation}
\label{sec:loewner-interpolation}

In this subsection, we briefly recall the Loewner framework for matrix interpolation according to \cite{mayo2007}. We refer the reader to the handbooks \cite{antoulas2017,gosea2022} for more in-depth details on various applications and extensions.

In principle, there are three different formulations of the Loewner framework depending on the input data, namely \emph{scalar} data, \emph{matrix} data, and so-called \emph{tangential} data. As the tangential formulation generalizes the other formulations it is most commonly used throughout the literature. For the sake of simplicity, we will use the Loewner framework for matrix data in this work.

Given transfer function samples in the form of pairs of scalars and matrices
\begin{equation*}
\left\{\bigl(p_i,G(p_i)\bigr)~|~p_i\in\Pi,\,G(p_i)\in\mathbb{R}^{k_1\times k_2},\,i\in\{1,\,\dots,\,n_p\right\},
\end{equation*}
where \(\Pi\subset P\) is the set of all parameter samples \(p_i\),
the Loewner framework aims to construct an interpolant \(\hat{G}(p)\) that fulfills the interpolation condition
\begin{equation}
\label{eq:loewner-interpolation}
\hat{G}(p_i)=G(p_i)\quad\text{for all}\quad i=1,\,\dots,\,n_p.
\end{equation}

The Loewner interpolant \(\hat{G}(p)\) is constructed by choosing a partitioning of the parameter samples (which is disjoint in our case)
\begin{equation}
\label{eq:partitioning}
    \{\pi_1,\,\dots,\,\pi_M\}\cup\{\phi_1,\,\dots,\,\phi_N\}=\Pi\subset \mathbb{R},
  \end{equation}
with \(N+M=n_p\), whereafter the (partitioned) transfer function samples are arranged in two matrices
\begin{equation}
\label{eq:loewner-quadruple-1}
\begin{aligned}
  \mathcal{V}&=
  \begin{bmatrix}
    G(\pi_1)^T&\cdots&G(\pi_M)^T
  \end{bmatrix}^T\in\mathbb{R}^{Mk_1\times k_2},\\
  \mathcal{W}&=
  \begin{bmatrix}
    G(\phi_1)&\cdots&G(\phi_N)
  \end{bmatrix}\in\mathbb{R}^{k_1\times Nk_2},
\end{aligned}
\end{equation}
and define the so-called \emph{Loewner matrix} \(\mathbb{L}\in\mathbb{R}^{Mk_1\times Nk_2}\) and \emph{shifted Loewner matrix} \(\mathbb{L}_s\in\mathbb{R}^{Mk_1\times Nk_2}\) as
\begin{equation}
\label{eq:loewner-quadruple-2}
\begin{aligned}
  \left\{\mathbb{L}\right\}_{i,j}&=\frac{G(\pi_i)-G(\phi_j)}{\pi_{i}-\phi_{j}}\in\mathbb{R}^{k_1\times k_2},\\
  \left\{\mathbb{L}_s\right\}_{i,j}&=\frac{\pi_{i}G(\pi_i)-\phi_{j}G(\phi_j)}{\pi_{i}-\phi_{j}}\in\mathbb{R}^{k_1\times k_2},
\end{aligned}
\end{equation}
respectively. The notations $\left\{\mathbb{L}\right\}_{i,j}$ and  $\left\{\mathbb{L}_s\right\}_{i,j}$ above refer to the $(i,j)$-th $k_1 \times k_2$ matrix block and not to the $(i,j)$-th scalar entry of each matrix, as it is conventional \cite{antoulas2017,mayo2007}.

Provided that
\begin{equation*}
\operatorname{rank}{p\mathbb{L}-\mathbb{L}_s}=\operatorname{rank}{\begin{bmatrix}\mathbb{L}&\mathbb{L}_s\end{bmatrix}}=\operatorname{rank}{\begin{bmatrix}\mathbb{L}\\\mathbb{L}_s\end{bmatrix}},
\end{equation*}
holds for all \(p\in\Pi\), we can compute the Singular Value Decompositions (SVDs) \cite{golub2013}
\begin{equation}
\label{eq:loewner-svd}
  X\Sigma_1\tilde{Y}^T =
  \begin{bmatrix}
    \mathbb{L}&\mathbb{L}_s
  \end{bmatrix},\quad\tilde{X}\Sigma_2Y^T=
  \begin{bmatrix}
    \mathbb{L}\\\mathbb{L}_s
  \end{bmatrix},
\end{equation}
and choose $r=\operatorname{rank}{\begin{bmatrix}\mathbb{L}&\mathbb{L}_s\end{bmatrix}}$. The truncation rank $r$ is then used to truncate $X$ and $Y$ by keeping the first $r$ columns of each matrix \(X=[X_r~\,*\,]\) and \(Y=[Y_r~\,*\,]\) with \(X_r\in\mathbb{R}^{Mk_1\times r}\) and \(Y_r\in\mathbb{R}^{Nk_2\times r}\) in order to construct a truncated transfer function
\begin{equation}
\label{eq:loewner-realization}
  \hat{G}(p)=\mathcal{W} Y_r(X_r^T(\mathbb{L}_s-p\mathbb{L})Y_r)^{-1}X_r^T\mathcal{V}.
\end{equation}
From \cite[Lemma 5.2]{mayo2007}, it follows that \(\hat{G}(p)\) tangentially interpolates the data, i.e. \eqref{eq:loewner-interpolation} holds.

So far, \(\hat{G}(p)\) according to \eqref{eq:loewner-realization} only interpolates the parametric samples as a function that is independent of the frequency variable $s$. In the following section, it is shown how the Loewner framework and the LFT can be conjointly used to construct a global bivariate transfer function and a corresponding state-space representation.

\section{Snapshot interpolation}
\label{sec:method}
The proposed matrix interpolation scheme based on the classical Loewner framework and LFT will be introduced in this section.

Assume that the input data is given in the form of local non-parametric models as in \eqref{eq:snapshots}. As mentioned in earlier, input data in descriptor form can be transformed by following the steps outlined in \cite{kunkel2006}.
Further, assume that all snapshots share a common state-space for all parameter values \(p_i\in\Pi\).  The latter means that the local snapshot realizations were not independently reduced or that a reprojection step onto suitable generalized state coordinates as described in \cite{geuss2013} was performed after the model reduction of individual snapshots. This preprocessing step is particularly important to avoid numerical issues when using different model reduction methods and different bases for the individual snapshots. This is discussed in more detail in \cite{geuss2013,yue2019}.

\begin{figure}
\centering
\begin{center}
\begin{tikzpicture}
\begin{axis}[
    width=.7\textwidth,
    view={45}{15},
    axis line style = {draw=none},
    tick style = {draw=none},
    ticks=none,
    axis equal image,
    xmin=0, xmax=7.5,
    ymin=-2*\gap, ymax=\n+\m+2*\gap,
    zmin=-2*\gap-.8-\gap, zmax=\n+\p+2*\gap,
    axis background/.style={fill=white},
    ]
\draw[->] (0,.5*\n+.5*\m,-.85) -- (7.5,.5*\n+.5*\m,-.85) node[below,pos=0.5] {\(p\)};
\snapshot{0}{1}
\snapshot{1.5}{2}
\snapshot{3}{3}
\draw[decorate sep={.7mm}{1.5mm},fill=white](3,.5*\n+.5*\m,.5*\n+.5*\p) -- (7.5,.5*\n+.5*\m,.5*\n+.5*\p);
\snapshot{7.5}{n}
\end{axis}
\end{tikzpicture}
\end{center}
\caption{Arrangement of the snapshot data for interpolation with the univariate Loewner framework.}
\label{fig:snapshots}
\end{figure}
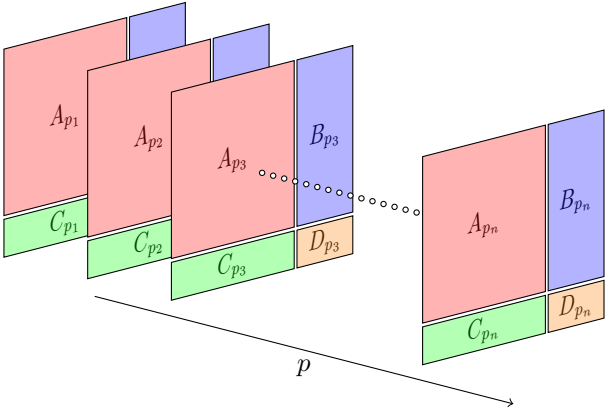

As insinuated in \cref{sec:lft}, we propose to arrange the snapshot realizations into larger matrices
\begin{equation}
\label{eq:snapshots-mat}
  G(p_i)=
  \begin{bmatrix}
    A(p_i) & B(p_i) \\ C(p_i) & D(p_i)
  \end{bmatrix}\in\mathbb{R}^{(n+n_{\mathrm{o}})\times(n+n_{\mathrm{i}})},
\end{equation}
which are used as input data \(\left\{(p_i,G(p_i))~|~i=1,\,\dots,\,n_p\right\}\) in the Loewner framework. The arrangement of the snapshot data is schematically visualized in \cref{fig:snapshots}.

By applying the Loewner framework and following the steps outlined in \cref{sec:loewner-interpolation} with partitioning \eqref{eq:partitioning}, a truncation rank \(r\) is chosen
to obtain an interpolant \(\hat{G}(p)\) via \eqref{eq:loewner-realization} that fulfills the interpolatory conditions \eqref{eq:loewner-interpolation}.
Subsequently, a parametric transfer function is constructed via an LFT (and recalling \cref{rem:schur}) as
\begin{equation}
\label{eq:ptf-lft}
    \hat{H}(s,p)=\mathcal{F}_{u}\left(\hat{G}(p),s^{-1}I_n\right).
\end{equation}
It should be noted that the interpolation is only performed in the parameter $p$ using univariate Loewner matrices. The bivariate transfer function is then obtained by applying the LFT.
To provide an explicit formula for \eqref{eq:ptf-lft}, we now introduce intermediate quantities that are useful in the remainder of this work. To this end, the definition of \(\hat{G}(p)\) is expanded according to \eqref{eq:loewner-realization} into a block form
\begin{align}
\label{eq:Gp_exp}
\hat{G}(p)&=\underbrace{\mathcal{W} Y_r}_{=:\left[\substack{\mathcal{Y}\\\mathcal{C}}\right]}
\underbrace{(X_r^T(\mathbb{L}_s-p\mathbb{L}) Y_r)^{-1}}_{=:\mathcal{K}(p)^{-1}}\underbrace{X_r^T\mathcal{V}}_{=:[\mathcal{X}\;\mathcal{B}]}\nonumber\\
&=\begin{bmatrix}\mathcal{Y}\mathcal{K}(p)^{-1}\mathcal{X}~&\mathcal{Y}\mathcal{K}(p)^{-1}\mathcal{B}\\\mathcal{C}\mathcal{K}(p)^{-1}\mathcal{X}~&\mathcal{C}\mathcal{K}(p)^{-1}\mathcal{B}\end{bmatrix},
\end{align}
where 
\begin{equation}
\label{eq:Kp}
\mathcal{K}(p)=X_r^T(\mathbb{L}_s-p\mathbb{L})Y_r\in\mathbb{R}^{r\times r}
\end{equation}
and the matrices \(\mathcal{Y}\in\mathbb{R}^{n\times r}\), \(\mathcal{C}\in\mathbb{R}^{n_{\mathrm{o}}\times r}\), \(\mathcal{X}\in\mathbb{R}^{r\times n}\) and \(\mathcal{B}\in\mathbb{R}^{r\times n_{\mathrm{i}}}\) are obtained by conformingly splitting \(\mathcal{W} Y_r\in\mathbb{R}^{(n+n_{\mathrm{o}})\times r}\) and \(X_r^T\mathcal{V}\in\mathbb{R}^{r\times(n+n_{\mathrm{i}})}\), respectively, i.e.
\begin{equation}
\label{eq:data}
\begin{aligned}
\mathcal{X}^T\!&=\begin{bmatrix}A(\pi_1)^T&\!C(\pi_1)^T&\cdots&A(\pi_M)^T&\!C(\pi_M)^T\end{bmatrix}X_r,\\
\mathcal{B}^T\!&=\begin{bmatrix}B(\pi_1)^T&\!\!D(\pi_1)^T&\cdots&B(\pi_M)^T&\!D(\pi_M)^T\end{bmatrix}X_r,\\
\mathcal{Y}&=\begin{bmatrix}A(\phi_1)&B(\phi_1)&\cdots&A(\phi_N)&B(\phi_N)\end{bmatrix}Y_r,\\
\mathcal{C}&=\begin{bmatrix}C(\phi_1)&\!D(\phi_1)&\cdots&C(\phi_N)&D(\phi_N)\end{bmatrix}Y_r.
\end{aligned}
\end{equation}
The realization of \(\hat{H}(s,p)\) according to \eqref{eq:ptf-lft} only requires $\hat{G}(p)$ and the LFT which is realization independent and does not depend on snapshot data. In other words, the realization of \(\hat{H}(s,p)\) can be solely constructed according to \eqref{eq:Gp_exp} from the snapshot data \eqref{eq:snapshots} and the Loewner projection matrices \(X_r\in\mathbb{R}^{M(n+n_{\mathrm{o}})\times r}\) and \(Y_r\in\mathbb{R}^{N(n+n_{\mathrm{i}})\times r}\). The following theorem shows how this can be done explicitly and provides interpolation conditions for the bivariate transfer function.
\begin{theorem}
\label{thm:general-real}
Let \(H(s,p)\) denote the parametric transfer function realized as \eqref{eq:ptf}. Furthermore, let \(\Pi\subset P\) denote a set of parameter samples and let \(G(p_i)\) denote the snapshot matrix at the parameter \(p_i\in\Pi\) according to \eqref{eq:snapshots-mat}.

If the rational function \(\hat{G}(p)\) as in \eqref{eq:loewner-realization} obtained by the Loewner framework from the input data \(\{(p_i,G(p_i))\}_{i=1}^{n_p}\) satisfies the interpolation condition in \eqref{eq:loewner-interpolation}, the parametric transfer function given by
\begin{equation}
\label{eq:precise-formula}
\hat{H}(s,p)=\hat{C}(p)(sI_n-\hat{A}(p))^{-1}\hat{B}(p)+\hat{D}(p)
\end{equation}
with
\begin{align*}
    \hat{A}(p)&=\mathcal{Y}\mathcal{K}(p)^{-1}\mathcal{X},&&\hat{B}(p)=\mathcal{Y}\mathcal{K}(p)^{-1}\mathcal{B},\\
    \hat{C}(p)&=\mathcal{C}\mathcal{K}(p)^{-1}\mathcal{X},&&\hat{D}(p)=\mathcal{C}\mathcal{K}(p)^{-1}\mathcal{B},
\end{align*}
satisfies
    \begin{equation}
    \label{eq:ptf-interpolation}
        \hat{H}(s,p_i)=H(s,p_i)
    \end{equation}
    for \((s,p_i)\in\mathbb{C}\times\Pi\).
\end{theorem}
\begin{proof}
Incorporating the reformulated expression \eqref{eq:Gp_exp} into the LFT in \eqref{eq:ptf-lft} yields
\begin{align*}
&\hat{H}(s,p)=\mathcal{F}_{u}\left(\begin{bmatrix}\mathcal{Y}\mathcal{K}(p)^{-1}\mathcal{X}&\mathcal{Y}\mathcal{K}(p)^{-1}\mathcal{B}\\\mathcal{C}\mathcal{K}(p)^{-1}\mathcal{X}&\mathcal{C}\mathcal{K}(p)^{-1}\mathcal{B}\end{bmatrix},s^{-1}I_n\right)\nonumber\\
&=\mathcal{C}\mathcal{K}(p)^{-1}\mathcal{X}(sI_n\!-\!\mathcal{Y}\mathcal{K}(p)^{-1}\mathcal{X})^{-1}\mathcal{Y}\mathcal{K}(p)^{-1}\mathcal{B}+\mathcal{C}\mathcal{K}(p)^{-1}\mathcal{B}.
\end{align*}
The resulting matrices in \eqref{eq:precise-formula} are then trivial.
In addition, recalling \cref{rem:schur}, the interpolation condition \eqref{eq:ptf-interpolation} immediately follows from \eqref{eq:loewner-interpolation} and \eqref{eq:ptf-lft}.
\end{proof}
\cref{thm:general-real} provides a formula for the bivariate transfer function in a general form that depends on frequency and on the parameter. Under certain conditions, a more compact form of \eqref{eq:precise-formula} can be obtained.

\begin{theorem}
\label{thm:compact-real}
 Let $\hat{H}(s,p)$ denote the parametric transfer function in \cref{thm:general-real}. If \(\mathcal{K}(p)\) and \((\mathcal{K}(p)-s^{-1}\mathcal{X}\mathcal{Y})\) are invertible, then
 \(s\neq0\) and $\hat{H}(s,p)$ can be written as
 \begin{equation}
 \label{eq:compact-formula}
     \hat{H}(s,p)= \mathcal{C}(p\mathcal{E}-s^{-1}\mathcal{X}\mathcal{Y}-\mathcal{A})^{-1}\mathcal{B},
 \end{equation}
 where 
\begin{align*}
\mathcal{E}&=-X_r^T\mathbb{L} Y_r\in\mathbb{R}^{r\times r},\\
\mathcal{A}&=-X_r^T\mathbb{L}_s Y_r\in\mathbb{R}^{r\times r},\\\mathcal{X}\mathcal{Y}&=X_r^T\mathbb{T} Y_r\in\mathbb{R}^{r\times r},
\end{align*}
with
\begin{equation}
\label{eq:m}
\{\mathbb{T}\}_{i,j}=\begin{bmatrix}A(\pi_i)A(\phi_j)&A(\pi_i)B(\phi_j)\\C(\pi_i)A(\phi_j)&C(\pi_i)B(\phi_j)\end{bmatrix}.
\end{equation}
\end{theorem}
\begin{proof}
The invertibility of \((\mathcal{K}(p)-s^{-1}\mathcal{X}\mathcal{Y})\) implies \(s\neq0\) by contradictory argument. If \(\mathcal{K}(p)\) is also invertible, we can invoke the Sherman-Woodbury-Morrison matrix identity \cite[Sec. 2.1.4.]{golub2013}:
\begin{equation*}
        (S+KL^*)^{-1}=S^{-1}-S^{-1}K(I_n+L^*S^{-1}K)^{-1}L^*S^{-1}
\end{equation*}
by choosing \(S=\mathcal{K}(p)\), \(K=-s^{-1}\mathcal{X}\) and \(L=\mathcal{Y}^*\) in \eqref{eq:precise-formula} to obtain
\begin{equation*}
    \hat{H}(s,p)=\mathcal{C}\left(\mathcal{K}(p)-s^{-1}\mathcal{X}\mathcal{Y}\right)^{-1}\mathcal{B}.
\end{equation*}
Then, we expand \(\mathcal{K}(p)=X_r^T(\mathbb{L}_s-p\mathbb{L})Y_r\) and set \(\mathcal{A}=-X_r^T\mathbb{L}_s Y_r\), \(\mathcal{E}=-X_r^T\mathbb{L} Y_r\) with \(\mathbb{L}\) and \(\mathbb{L}_s\) given by \eqref{eq:loewner-quadruple-2} to get the final expression,
\begin{equation*}
  \hat{H}(s,p)=\mathcal{C}\left(p\mathcal{E}-s^{-1}\mathcal{X}\mathcal{Y}-\mathcal{A}\right)^{-1}\mathcal{B}.
\end{equation*}
The identity
\begin{equation*}
    \mathcal{X}\mathcal{Y}=X_r^T\begin{bmatrix}
        A(\pi_1)\\C(\pi_1)\\\vdots\\A(\pi_M)\\C(\pi_M)
    \end{bmatrix} 
    \begin{bmatrix}
        A(\phi_1)^T\\B(\phi_1)^T\\\vdots\\A(\phi_N)^T\\B(\phi_N)^T
    \end{bmatrix}^TY_r=X_r^T\mathbb{T} Y_r
\end{equation*}
with \(\mathbb{T}\) as in \eqref{eq:m} can be verified by checking \eqref{eq:Gp_exp}.
\end{proof}
\begin{remark}
The bivariate transfer function can be rewritten in order to avoid the inverse of $s$ such that
    \begin{align*}
    \hat{H}(s,p)&=\mathcal{C}\left(p\mathcal{E}-s^{-1}\mathcal{X}\mathcal{Y}-\mathcal{A}\right)^{-1}\mathcal{B}\\
    &=s\mathcal{C}(sp\mathcal{E}-s\mathcal{A}-\mathcal{X}\mathcal{Y})^{-1}\mathcal{B}.
\end{align*}
\end{remark}
\begin{remark}
Note that if $s=0$, by plugging in \eqref{eq:Gp_exp} into \eqref{eq:schur}, it follows that
\begin{equation}
\begin{aligned}
    \label{eq:zero-formula}
    \hat{H}(0,p)=\,&\mathcal{C}\mathcal{K}(p)^{-1}\mathcal{B}\\&-\mathcal{C}\mathcal{K}(p)^{-1}\mathcal{X}(\mathcal{Y}\mathcal{K}(p)^{-1}\mathcal{X})^{-1}\mathcal{Y}\mathcal{K}(p)^{-1}\mathcal{B}
\end{aligned}
\end{equation}
Along the lines of \cite{beattie2009}, we define \(\mathcal{P}(p)=\mathcal{X}(\mathcal{Y}\mathcal{K}(p)^{-1}\mathcal{X})^{-1}\allowbreak\mathcal{Y}\mathcal{K}(p)^{-1}\) and observe that \(\mathcal{P}(p)\mathcal{P}(p)=\mathcal{P}(p)\) which means that \(\mathcal{P}(p)\) is an oblique projection onto \(\operatorname{ran}(\mathcal{X})\). We can arrive at a concise expression by substituting \(\mathcal{P}(p)\) in \eqref{eq:zero-formula}:
\begin{align}
    \hat{H}(0,p)&=\mathcal{C}\mathcal{K}(p)^{-1}\mathcal{B}-\mathcal{C}\mathcal{K}(p)^{-1}\mathcal{P}(p)\mathcal{B}\nonumber\\
    &=\mathcal{C}\mathcal{K}(p)^{-1}(I-\mathcal{P}(p))\mathcal{B}\label{eq:formula-schur}.
\end{align}
\end{remark}
 The steps to obtain the parametric realization of \(\hat{H}(s,p)\) based on snapshot data of the form \eqref{eq:snapshots} are summarized in \cref{alg:1}. One notable advantage of this method is that the coefficient matrices of the final interpolated parametric transfer function can be directly constructed based on snapshot data without additional computations in the online phase. Furthermore, the matrices are truncated based on the rank of the Loewner and the shifted Loewner matrices which further compresses the final model. This is further discussed in \cref{sec:rank}. The numerical properties of the formulas offered by \cref{thm:general-real} and \cref{thm:compact-real} are discussed in \cref{sec:numerics}.

\begin{algorithm2e}
\caption{Parametric interpolation}\label{alg:1}
\Input{Set of parameter samples \(\Pi\), snapshot realizations \(G(p_i)\) as in \eqref{eq:snapshots-mat} for every \(p_i\in\Pi\), truncation tolerance \(\varepsilon>0\).}
\Output{\((\mathcal{E},\mathcal{A},\mathcal{B},\mathcal{C},\mathcal{X},\mathcal{Y})\).}
Partition the parameter samples into two disjoint sets as in \eqref{eq:partitioning}.\\
Construct the Loewner and shifted Loewner matrices according to \eqref{eq:loewner-quadruple-2}.\\
Obtain the left and right subspaces \(X,Y\) by computing the SVDs \eqref{eq:loewner-svd}.\\
Choose a truncation rank \(r\) with respect to the tolerance \(\varepsilon>0\) according to
\begin{equation}
\label{eq:loewner-tol}
    r=\min\left\{~k~\left|\left(\varepsilon\leq\frac{\sum_{i=1}^{k-1}|\sigma_i|^2}{\sum_{i=1}^n|\sigma_i|^2}\right)^{1/2}\right.\right\}
\end{equation}
and compute the truncated Loewner singular vectors \(X_r\) and \(Y_r\) from \eqref{eq:loewner-svd}.\\
Construct the matrices \(\mathcal{V}\) and \(\mathcal{W}\) from the data as in \eqref{eq:loewner-quadruple-1} and compute and partition
\begin{align*}
X_r^T\mathcal{V}=\begin{bmatrix}\mathcal{X}&\mathcal{B}\end{bmatrix}~\text{and}~\mathcal{W} Y_r=\begin{bmatrix}\mathcal{Y}\\\mathcal{C}\end{bmatrix}.
\end{align*}\\
Compute the remaining matrices via
\begin{align*}
\mathcal{E}&=-X_r^T\mathbb{L} Y_r,&&\mathcal{A}=-X_r^T\mathbb{L}_s Y_r.
\end{align*}
\end{algorithm2e}

\section{Rank bounds on \texorpdfstring{\(\mathbb{L}\)}{L} and \texorpdfstring{\(\mathbb{L}_s\)}{Ls}}
\label{sec:rank}
The ranks of the Loewner and shifted Loewner matrices are indicators of the truncation order. For the multivariate Loewner approach, introduced in \cite{antoulas2012}, and further developed in \cite{ionita2014}, the rank of the multivariate Loewner matrix scales with the dimension of the rational interpolant, but also with the dimension of the input data. This property is inherently different than for the univariate Loewner approach, in which the rank of the Loewner matrix indicates the dimension of the minimal rational interpolant (assuming that enough data is provided).

Although it is a challenging task to know a priori what the exact ranks of the Loewner and shifted Loewner matrices would be, one can infer a rank bound from the specific kind of parametric structure of the transfer function. In what follows, we restrict our discussion to polynomial dependence on \(p\) as we can break down the Loewner matrices in an affine way in terms of the coefficients which allows us to obtain rank bounds. In this section, we distinguish two cases: an affine dependence of $G(p)$ on \(p\) and a higher-order polynomial dependence on \(p\).

\subsection{Affine dependence on \texorpdfstring{$p$}{p}}
An affine parameter dependence of the realization matrices often appears in practical scenarios, see \cite{benner2015,ionita2014}, and also the two numerical examples treated in \cref{sec:numerics}. In this case, the rational matrix $G(p)$ is considered to be a polynomial matrix that is linear in $p$ and thus can be written as
\begin{equation}
\label{eq:G-affine}
G(p)= \Gamma_0 + p \Gamma_1,
\end{equation}
where $\Gamma_0,\Gamma_1\in\mathbb{R}^{(n+n_{\mathrm{o}})\times(n+n_{\mathrm{i}})}$ are matrix coefficients:
Then, it follows that the $(i,j)$-th blocks of the Loewner and shifted Loewner matrices $\mathbb{L}, \mathbb{L}_s \in \mathbb{R}^{M(n+n_{\mathrm{o}})\times N(n+n_{\mathrm{i}})}$ can be expressed as
\begin{equation}
\label{eq:LandLs}
\begin{aligned}
    \left\{\mathbb{L}\right\}_{i,j} &= \frac{G(\pi_i)-G(\phi_j)}{\pi_i-\phi_j} = \frac{\Gamma_1 \pi_i - \Gamma_1 \phi_j}{\pi_i - \phi_j} = \Gamma_1, \\
        \left\{\mathbb{L}_s\right\}_{i,j} &= \frac{\pi_i G(\pi_i)-\phi_j G(\phi_j)}{\pi_i-\phi_j}\\&= \frac{\pi_i(\Gamma_0+\Gamma_1 \pi_i) - \phi_j(\Gamma_0+\Gamma_1 \phi_j)}{\pi_i - \phi_j} \\&= \Gamma_0+\Gamma_1(\pi_i+\phi_j),
    \end{aligned}
\end{equation}
where $ \left\{\mathbb{L}\right\}_{i,j}$ and $ \left\{\mathbb{L}_s\right\}_{i,j}$ represent $(n+n_{\mathrm{o}}) \times (n+n_{\mathrm{i}})$ blocks of the Loewner matrices, associated with the sample pair $(\pi_i,\phi_j)$. Then, based on the above formulations, the Loewner matrices can be written in a compact Kronecker format as stated in \cref{Lemma:LLs}.
\begin{lemma}
\label{Lemma:LLs}
If \(G(p)\) as in \eqref{eq:G-affine} depends affinely on \(p\), the Loewner matrices \(\mathbb{L}\) and \(\mathbb{L}_s\) can be explicitly expressed as
\begin{align*}
    \mathbb{L} = \mathbb{1}_{M,N}  \otimes    \Gamma_1, &&\mathbb{L}_s =  \mathbb{1}_{M,N}  \otimes    \Gamma_0 + (\mathbb{1}_{M,N}  \otimes    \Gamma_1) \odot \Xi_2,
    \end{align*}
    where $\mathbb{1}_{M,N}$ is a $M \times N$ matrix of ones, $\otimes$ and $\odot$ denote the Kronecker and Hadamard products \cite{golub2013}, respectively, and $\Xi_2 \in \mathbb{R}^{M(n+n_{\mathrm{o}})\times N(n+n_{\mathrm{i}})}$ is given by $ \left\{\Xi_2 \right\}_{i,j}= (\pi_i+\phi_j) \mathbb{1}_{n+n_{\mathrm{o}},n+n_{\mathrm{i}}}$ and can be written as
    \begin{equation}
    \label{eq:xi2}
        \Xi_2 = \pi \otimes \mathbb{1}_{(n+n_{\mathrm{o}}) \times N(n+n_{\mathrm{i}})} + \phi^T \otimes \mathbb{1}_{M(n+n_{\mathrm{o}}) \times (n+n_{\mathrm{i}})},
    \end{equation}
    where $\pi = \begin{bmatrix}
        \pi_1 & \pi_2 & \cdots & \pi_M
    \end{bmatrix}^T$,  $\phi = \begin{bmatrix}
        \phi_1 & \phi_2 & \cdots & \phi_N
    \end{bmatrix}^T$.
\end{lemma}
\begin{proof}
The result can be directly deduced from the formulas stated in \eqref{eq:LandLs}.
\end{proof}
\cref{Lemma:LLs} is now used to derive rank bounds on the Loewner and shifted Loewner matrices in the following theorem.
\begin{theorem}
\label{thm:ranks}
Let \(G(p)\) as in \eqref{eq:G-affine} have an affine parametric dependence. Then,
    \begin{align*}
        \operatorname{rank}(\mathbb{L})&= \operatorname{rank}(\Gamma_1), \\  \operatorname{rank}(\mathbb{L}_s) &\leq \operatorname{rank}(\Gamma_0) + \operatorname{rank}(\Gamma_1) \operatorname{rank}(\Xi_2),
    \end{align*}
    where \(\Xi_2\) is given by \eqref{eq:xi2}.
\end{theorem}
\begin{proof}
The result can be directly deduced from \cref{Lemma:LLs} by using the rank inequalities
\begin{equation}
\label{eq:rank-ineqs}
\begin{split}
    \operatorname{rank}(X+Y)&\leq\operatorname{rank}(X)+\operatorname{rank}(Y)\\
    \operatorname{rank}(X\odot Y)&\leq\operatorname{rank}(X)\operatorname{rank}(Y)
\end{split}
\end{equation}
together with \(\operatorname{rank}{\mathbb{1}_{M,N}\otimes  X}= \operatorname{rank}(X)\). A proof for the rank inequality involving the Hadamard product is given in \cite{ballantine1968}.
\end{proof}
\cref{thm:ranks} can be used to obtain a priori rank bounds on the Loewner and shifted Loewner matrices based on the relation of the transfer function to the parameter in the case of affine dependence. In what follows we extend these results for general higher order polynomial dependence. As for the classical parametric Loewner approach \cite{antoulas2012,ionita2014}, for which the rank of the multi-dimensional Loewner matrices is shown to depend not only on the dimension of the minimal interpolant, a dependence on the chosen data is also noticeable. In what follows, this is encoded in the appearance of matrices $\Xi_k$, which are written solely in terms of the interpolation sample points.

\subsection{Higher-order polynomial dependence}
More generally, \(G(p)\) can also be a polynomial matrix of degree \(h\in\mathbb{N}\)
\begin{equation}
\label{eq:G-polynomial}
    G(p)=\sum_{k=0}^h p^k\Gamma_k
\end{equation}
with matrix coefficients \(\Gamma_k\in\mathbb{R}^{(n+n_{\mathrm{o}})\times(n+n_{\mathrm{i}})}\). The rank bounds presented in the previous subsection are generalized for the polynomial case.

\begin{theorem}
\label{thm:ranks_gen}
Let \(G(p)\) as in \eqref{eq:G-polynomial} have a polynomial dependence on \(p\). The following inequalities hold:
    \begin{align*}
        \operatorname{rank}(\mathbb{L}) &\leq \sum_{k=0}^h \operatorname{rank}(\Gamma_k) \operatorname{rank}(\Xi_k), \\
        \operatorname{rank}(\mathbb{L}_s) &\leq \sum_{k=0}^h \operatorname{rank}(\Gamma_k) \operatorname{rank}(\Xi_{k+1}),
    \end{align*}
    where \(\Xi_k\in\mathbb{R}^{M(n+n_{\mathrm{o}})\times N(n+n_{\mathrm{i}})}\) is given by
    \begin{equation*}
        \{\Xi_k\}_{i,j}=\left(\frac{\pi_i^k-\phi_j^k}{\pi_i - \phi_j} \right) \mathbb{1}_{n+n_{\mathrm{o}},n+n_{\mathrm{i}}}.
    \end{equation*}
\end{theorem}

\begin{proof}
By using the following identities 
\begin{align*}
    \left\{\mathbb{L}\right\}_{i,j} &=  \frac{\sum_{k=0}^h \pi_i^k\Gamma_k - \sum_{k=0}^h \phi_j^k\Gamma_k}{\pi_i - \phi_j} \\&= \sum_{k=0}^h \frac{\pi_i^k-\phi_j^k}{\pi_i - \phi_j} \Gamma_k, \\
        \left\{\mathbb{L}_s\right\}_{i,j} &= \frac{\sum_{k=0}^h \pi_i^{k+1}\Gamma_k - \sum_{k=0}^h \phi_j^{k+1}\Gamma_k}{\pi_i - \phi_j} \\&= \sum_{k=0}^h \frac{\pi_i^{k+1}-\phi_j^{k+1}}{\pi_i - \phi_j} \Gamma_k,
\end{align*}
together with the classical rank inequalities \eqref{eq:rank-ineqs}, the two inequalities above directly follow.
\end{proof}
\begin{remark}
By definition, $\Xi_0$ is an all-zero matrix and $\Xi_1$ is an all-one matrix, while $\Xi_2$ is precisely defined as in \eqref{eq:xi2}.
The inequalities in \cref{thm:ranks_gen} can be viewed as a direct generalization
of those in \cref{thm:ranks}. To see this, one can rewrite the former by using that \(\operatorname{rank}(\Xi_0) = 0\) and \(\operatorname{rank}(\Xi_1) = 1\) as:
\begin{align*}
    \operatorname{rank}(\mathbb{L}) &\leq \operatorname{rank}(\Gamma_1) + \sum_{k=2}^h \operatorname{rank}(\Gamma_k) \operatorname{rank}(\Xi_k), \\
    \operatorname{rank}(\mathbb{L}_s) &\leq \operatorname{rank}(\Gamma_0)  + \sum_{k=1}^h \operatorname{rank}(\Gamma_k) \operatorname{rank}(\Xi_{k+1}),
\end{align*}
and setting $h = 1$. \cref{thm:ranks} is stricter in the sense that the bound on \(\mathbb{L}\) becomes sharp.
\end{remark}

\subsection{Case Studies}
In this subsection, we present three case studies to validate the results of \cref{thm:ranks} and \cref{thm:ranks_gen}. Some of the examples introduced here will be revisited for numerical studies in \cref{sec:numerics}.
\begin{example}[Toy system]
\label{ex:toy}
Consider the following realization
\begin{align*}
A(p)=\begin{bmatrix*}[r]-2&p&0\\-p&-1&0\\0&0&-1\end{bmatrix*},&&
B(p)=\begin{bmatrix}
    1\\0\\1
\end{bmatrix},&&\begin{aligned}
C(p)&=B(p)^T,\\
D(p)&=0.
\end{aligned}
\end{align*}
Hence, for the block matrix $G(p)$ assembled like \eqref{eq:snapshots-mat}, we can identify the following polynomial matrix coefficients:
\begin{align*}
    \Gamma_0\!&=\!\scalebox{.8}{\(\begin{bmatrix*}[r] -2 & 0 & 0 & 1\\ 0 & -1 & 0 & 0\\ 0 & 0 & -1 & 1\\ 1 & 0 & 1 & \phantom{-}0 \end{bmatrix*}\)},
    &&\Gamma_1\!=\!\scalebox{.8}{\(\begin{bmatrix*}[r] 0 & \phantom{-}1 & \phantom{-}0 & \phantom{-}0\\ -1 & 0 & 0 & 0\\ 0 & 0 & 0 & 0\\ 0 & 0 & 0 & 0 \end{bmatrix*}\)}.
\end{align*}
In this case, $\operatorname{rank}(\Gamma_1)=2$ and $\operatorname{rank}(\Gamma_0)=4$. Next, choose
the left points \(\pi_i\in\{0.5,1.5\}\) and the right points \(\phi_j\in\{2,4\}\)and evaluate $G(p)$ at these points. It follows that $\operatorname{rank}(\mathbb{L})=2$, while $\operatorname{rank}(\mathbb{L}_s)=6$. Hence, it can be verified that the conditions in \cref{thm:ranks} are valid:
\begin{align*}
    \operatorname{rank}(\mathbb{L}) &= \operatorname{rank}(\Gamma_1) = 2,\\
    6 = \operatorname{rank}(\mathbb{L}_s) &< \operatorname{rank}(\Gamma_0) + \operatorname{rank}(\Gamma_1)\operatorname{rank}(\Xi_2)  = 8.
\end{align*}
\end{example}

\begin{example}[Modified toy system]
We modify the $(3,3)$ entry of the matrix $A(p)$, by replacing $-1$ with $-p$.

In this case, the $\operatorname{rank}{\Gamma_1}$ increases from $2$ to $3$, since its $(3,3)$ entry changes from $0$ to $-1$. The $(3,3)$ entry of matrix $\Gamma_0$ changes from $-1$ to $0$, but its rank remains $4$. Also, the rank of $\mathbb{L}$ jumps from $2$ to $3$, while the $\operatorname{rank}{\mathbb{L}_s}$ stays at 6.
So, it follows that
\begin{align*}
    \operatorname{rank}(\mathbb{L}) &= \operatorname{rank}(\Gamma_1) = 3,\\
    6 = \operatorname{rank}(\mathbb{L}_s) &< \operatorname{rank}(\Gamma_0) + \operatorname{rank}(\Gamma_1) \operatorname{rank}(\Xi_2)  = 10,
\end{align*}which again validates the identities of \cref{thm:ranks}.
\end{example}

\begin{example}[Polynomial system]
\label{ex:polynomial}
Consider the following example with polynomial dependence on $p$:
\begin{align*}
A(p)&=\begin{bmatrix}0.1p^2-2&p^3-p&0.2p^2\\-p^3&p^2-1&-0.5p\\-0.2p^2&-10p^3-0.5p&-1\end{bmatrix},&&
B(p)=\begin{bmatrix}
    1\\p\\1
\end{bmatrix},\\
C(p)&=\begin{bmatrix}
    1&0&1
\end{bmatrix},&&D(p)=0.
\end{align*}
We can directly identify the following polynomial matrix coefficients:
\begin{align*}
    \Gamma_0\!&=\!\scalebox{.8}{\(\begin{bmatrix*}[r]
    -2 & 0 & 0 & \phantom{.1}1\\
    0 & \phantom{.1}-1 & 0 & 0\\
    0 & 0 & \phantom{.}-1 & 1\\
    1 & 0 & 1 & 0 \end{bmatrix*}\)},
    &&\Gamma_1\!=\!\scalebox{.8}{\(\begin{bmatrix*}[r]
    0 & -1\phantom{.5} & \phantom{-}0\phantom{.5} & 0\\
    -1 & \phantom{-}0\phantom{.5} & -0.5 & 1 \\
    0 & -0.5 & \phantom{-}0\phantom{.5} & 0\\ 0
    & \phantom{-}0\phantom{.5} & \phantom{-}0\phantom{.5} & 0 \end{bmatrix*}\)}, \\
    \Gamma_2\!&=\!\scalebox{.8}{\(\begin{bmatrix*}[r]
    0.1 & \phantom{-}0 & \phantom{-}0.2 & \phantom{.1}0\\ 
    \phantom{-}0\phantom{.5} & 1 & \phantom{-}0\phantom{.5} & 0\\
    -0.2 & 0 & 0\phantom{.5} & 0\\ 
    0\phantom{.5} & 0 & 0\phantom{.5} & 0 \end{bmatrix*}\)},
    &&\Gamma_3\!=\!\scalebox{0.8}{\(\begin{bmatrix*}[r]
    0 & 1 & \phantom{-}0 & \phantom{-.}0\\
    -1 & 0 & 0 & 0 \\
    0 & -10 & 0 & 0\\
    0 & 0 & 0 & 0 \end{bmatrix*}\)}.
\end{align*}
It follows that \(\operatorname{rank}(\mathbb{L}_s)=4\), \(\operatorname{rank}(\mathbb{L}_s)=2\), \(\operatorname{rank}(\Gamma_2)=3\) and \(\operatorname{rank}(\Gamma_3)=2\). To ensure the full recovery of the original cubic model, more data needs to be used (compared to the previous two examples). This is indeed normal and to be expected, since the degree of the polynomial matrix increased from $1$ to $3$.

Hence, we choose the left points \(\pi_i\in\{0.5,1.5,2.5,3.5\}\) and the right points \(\phi_j\in\{2,4,6,8\}\). We form the Loewner matrices $\mathbb{L},\mathbb{L}_s\in\mathbb{R}^{16\times16}$ with the following properties: $\operatorname{rank}(\mathbb{L}) = 8$ and $\operatorname{rank}(\mathbb{L}_s) = 11$. Finally, we construct matrices $\Xi_k$ for $1\leq k \leq 4$ with $\operatorname{rank}(\mathbb{L}_s) = k$ for any $k$.

The next inequalities hold according to \cref{thm:ranks_gen}:
\begin{align*}
   8=\operatorname{rank}(\mathbb{L})&<\operatorname{rank}(\Gamma_1) &+\operatorname{rank}(\Gamma_2) \operatorname{rank}(\Xi_2)& \\
   &&+ \operatorname{rank}(\Gamma_3) \operatorname{rank}(\Xi_3) &= 14,\\
   11=\operatorname{rank}(\mathbb{L}_s)&<\operatorname{rank}(\Gamma_0) &+ \operatorname{rank}(\Gamma_1) \operatorname{rank}(\Xi_2)& \\
   &&+ \operatorname{rank}(\Gamma_2) \operatorname{rank}(\Xi_3)&\\
   && + \operatorname{rank}(\Gamma_3) \operatorname{rank}(\Xi_4) &= 25.
\end{align*}
\end{example}

\section{Numerical experiments}
\label{sec:numerics}
In this section, we present three examples of the form \eqref{eq:ptf}. Note that in the developed method and in these three test cases, $p$ is scalar.

The numerical experiments were performed on an Intel\textsuperscript{(R)} Core\textsuperscript{TM} i5-9400F with 16 GB RAM running on Linux Ubuntu 22.04.3 LTS. The algorithm described in \cref{sec:method} was implemented in the Python programming language (version 3.10.12) using the open source library pyMOR \cite{milk2016} (version 2023.2) and tested for different applications.
\begin{center}%
  \setlength{\fboxsep}{5pt}%
  \fbox{%
  \begin{minipage}{.92\linewidth}
    \textbf{Source code availability}\newline
    The source code of the implementations used for the numerical experiments is available at
    \begin{center}
      \href{https://zenodo.org/doi/10.5281/zenodo.11246111}%
        {\texttt{doi:10.5281/zenodo.11246111}}
    \end{center}
    under the MIT license, authored by Art J. R. Pelling.
  \end{minipage}}
\end{center}

In each numerical example, the parametric transfer function \eqref{eq:ptf} is sampled at different parametric sample points \(p_i\in\Pi\) and the corresponding snapshot data \eqref{eq:snapshots} is collected. This data is directly used to construct the interpolated parametric transfer function $\hat{H}(s,p)$ based on \cref{alg:1}.

Before moving to the numerical experiments, we consider a computation strategy for \(\hat{H}(s,p)\) by comparing the two formulas \eqref{eq:precise-formula} and \eqref{eq:compact-formula}. With the formula in \eqref{eq:precise-formula}, two matrix inversions, namely \(\mathcal{K}(p)^{-1}\in\mathbb{R}^{r\times r}\) and \((sI_n-\mathcal{Y}\mathcal{K}(p)^{-1}\mathcal{X})^{-1}\in\mathbb{C}^{n\times n}\), as well as several matrix multiplications are performed for every value of \((s,p)\). This can become costly, especially if the order \(n\) of the snapshots is large. In contrast, the compact formula in \eqref{eq:compact-formula} only requires a single inversion of \((\mathcal{K}(p)-s^{-1}\mathcal{X}\mathcal{Y})\in\mathbb{C}^{r\times r}\) and two matrix multiplications. However, the latter is obtained by invoking the Sherman-Morrison-Woodbury matrix identity \cite[Sec. 2.1.4.]{golub2013} which can cause numerical issues in practice.

Here, we face the dilemma of \emph{fast vs. feasible}. Although, the compact formula \eqref{eq:compact-formula} is computationally more efficient and by consequence a preferable choice, one has to ensure that the results are numerically accurate. To resolve this dilemma, we propose an ad hoc criterion for switching between the two formulas.

In order to obtain formula \eqref{eq:compact-formula} by \cref{thm:compact-real}, it is assumed that \((\mathcal{K}(p)-s^{-1}\mathcal{X}\mathcal{Y})\in\mathbb{C}^{r\times r}\) is invertible. In the following, we will denote this quantity by
\begin{equation*}
    Z(s,p)=\mathcal{K}(p)-s^{-1}\mathcal{X}\mathcal{Y}.
\end{equation*}
Note that \(Z(s,p)\) could be singular even for \(s\neq0\) depending on the system and data. From a computational standpoint, the accuracy of \eqref{eq:compact-formula} might deteriorate even for an invertible \(Z(s,p)\) when the matrix is close to being singular. Hence, the condition number for inversion \(\kappa(Z(s,p))=\lVert Z(s,p)^{-1}\rVert_2\lVert Z(s,p)\rVert_2\) is considered as an indicator of the numerical accuracy of \eqref{eq:compact-formula}. In case \(\kappa(Z(s,p))\) surpasses a preset tolerance $\varepsilon_\mathrm{cond}$, the precise formula in \eqref{eq:precise-formula} is used instead of the fast formula \eqref{eq:compact-formula}.

The proposed computation strategy is summarized in \cref{alg:Hsp}. Note that explicitly computing \(\kappa(Z(s,p))\) might be computationally more costly than simply using the precise formula in \eqref{eq:precise-formula} without switching. Therefore, we employ an estimator \(\tilde{\kappa}\) for the condition number in the \(1\)-norm \(\lVert Z(s,p)^{-1}\rVert_1\lVert Z(s,p)\rVert_1\) that estimates \(\lVert Z(s,p)^{-1}\rVert_1\) based on an LU decomposition by calling the LAPACK{} \cite{anderson1999} routine \texttt{zgecon} that can be evaluated at relatively low costs.

  \begin{algorithm2e}
\caption{Computation of \(\hat{H}(s,p)\)\label{alg:Hsp}}
\Input{\((s,p)\in\mathbb{C}\times P\),\\condition number tolerance \(\varepsilon_\mathrm{cond}\),\\\((\mathcal{E},\mathcal{A},\mathcal{B},\mathcal{C},\mathcal{X},\mathcal{Y})\) from \cref{alg:1}.}
\Output{Transfer function value \(\hat{H}(s,p)\)}
Form \(\mathcal{K}(p)=p\mathcal{E}-\mathcal{A}\)\\
Form \(Z=\mathcal{K}(p)-s^{-1}\mathcal{X}\mathcal{Y}\)\\
\uIf{\(\tilde{\kappa}(Z)\leq\varepsilon_\mathrm{cond}\)}{
Invert \(Z\in\mathbb{C}^{r\times r}\)\\
\Return \(\mathcal{C} Z^{-1}\mathcal{B}\)
}
\uElse{
Invert \(\mathcal{K}(p)\in\mathbb{R}^{r\times r}\)\\
Form \(Y=\mathcal{Y}\mathcal{K}(p)^{-1}\)\\
Invert \((sI_n-Y\mathcal{X})\in\mathbb{C}^{n\times n}\)\\
\Return \(\mathcal{C}(\mathcal{K}(p)^{-1}\mathcal{X}(sI_n-Y\mathcal{X})^{-1}Y+\mathcal{K}(p)^{-1})\mathcal{B}\)}
\end{algorithm2e}

\begin{example}[Toy system]
\label{ex:toy-num}
Consider the toy system introduced in \cref{ex:toy}. The realization matrices are sampled at four uniformly distributed parameters over the parameter interval \(P=[0,100]\), i.e. \(\Pi=\{0,\tfrac{100}{3},\tfrac{200}{3},100\}\). An alternating partitioning \cite{karachalios2021} is selected for the Loewner realization process. The spectra of the Loewner matrix, shifted Loewner matrix and Loewner pencil are depicted in \cref{fig:toy-svs} which reveal the intrinsic order of the model. A truncation rank \(r=6\) is chosen according to \eqref{eq:loewner-tol} with truncation tolerance \(\varepsilon=10^{-7}\).

\begin{figure}
\centering
\includegraphics{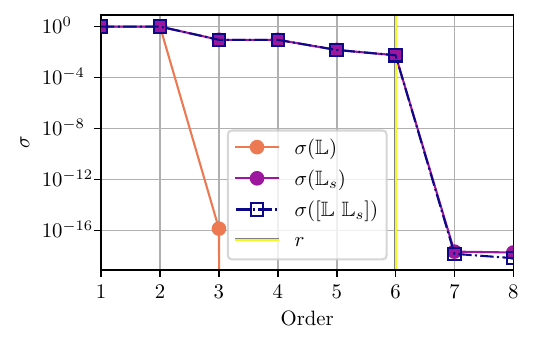}
\caption{Singular value decay of the Loewner matrices and pencil for the toy system (\cref{ex:toy-num}).}
\label{fig:toy-svs}
\end{figure}

After the computation of the parametric transfer function \(\hat{H}(s,p)\) based on \cref{alg:1}, the resulting transfer function is tested for ten equispaced parameter values in the same interval \(P\). The magnitude response of each test value is compared to the true frequency response $H(s,p)$ and depicted in \cref{fig:toy-mag}.

\begin{figure}
\centering
\includegraphics{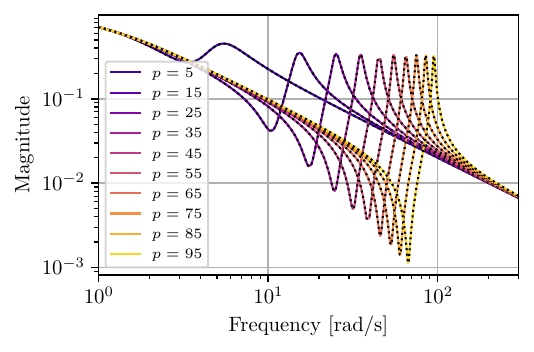}
\caption{Interpolated realizations \(\hat{H}(s,p)\) (solid colored lines) and true solutions \(H(s,p)\) (dotted black) for ten equispaced parameter values ranging from \(5\) to \(95\) for the toy system (\cref{ex:toy-num}).}
\label{fig:toy-mag}
\end{figure}

\begin{figure}
\centering
\includegraphics{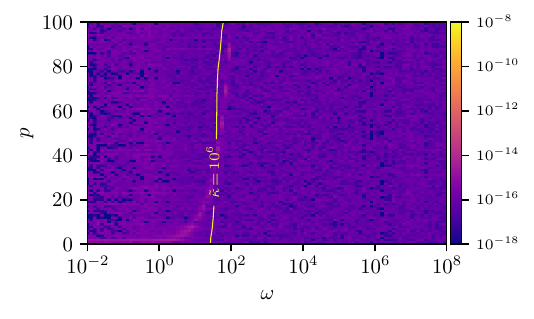}
\caption{Absolute error \eqref{eq:error} of the interpolated realizations over frequency and parameter for the toy system (\cref{ex:toy-num}).}
\label{fig:toy-2derr}
\end{figure}

A more systematic analysis of the interpolation quality is offered in \cref{fig:toy-2derr} which depicts the absolute error between the interpolated transfer function \(\hat{H}(s,p)\) and the true system
\begin{equation}
    \label{eq:error}
    \delta(\omega,p)=\lVert\hat{H}(\imath\omega,p)-H(\imath\omega,p)\rVert_2
\end{equation}
over frequency \(\omega\in\mathbb{R}\) and parameter \(p\in P\). A similar figure is found in \cite[Fig. 1.]{gosea2021} where the toy system was initially introduced. The yellow line in \cref{fig:toy-2derr} symbolizes the margin at which the switching of the different formulas described in \cref{alg:Hsp} is performed. The condition number in all examples increases towards lower frequencies, therefore, the precise formula \eqref{eq:precise-formula} is used to left of the line and the compact formula \eqref{eq:compact-formula} to the right. Values on the line correspond to values at which the condition number estimator is equal to the chosen tolerance, in this case \(\tilde{\kappa}(Z(s,p))=\varepsilon_\mathrm{cond}=10^{6}\). It was observed that the condition number estimator \(\tilde{\kappa}\) overestimates the true condition number \(\kappa(Z(s,p))\) which is about \(10^3\) on the yellow line in \cref{fig:toy-2derr}. Furthermore, the true condition number, as well as the estimator, seem to be mainly influenced by \(s\), growing rapidly towards \(s=0\) for the considered examples with an affine dependence on \(p\). For the polynomial system in \cref{ex:polynomial} the condition number varies with \(s\) and \(p\).
\cref{fig:toy-mag} and \cref{fig:toy-2derr} reveal that the interpolated parametric transfer function \(\hat{H}(s,p)\) is able to accurately represent the true system \(H(s,p)\) for parameter values \(p_i\) that are not part of the set of interpolation points \(\Pi\) as in \eqref{eq:partitioning}.
\end{example}

\begin{example}[Polynomial system]
\label{ex:polynomial-num}
In this example, we revisit the system with a polynomial dependency in the realization matrices from \cref{ex:polynomial} and sample the realization matrices at eight uniformly distributed parameters over the parameter interval \(P=[0,100]\). An alternating partitioning is selected for the Loewner realization process. According to \(\varepsilon=10^{-7}\), a truncation rank \(r=11\) is chosen.

\cref{fig:poly-2derr} shows the absolute error according to \eqref{eq:error} between interpolated and true magnitudes of the transfer function over a range of parameters and frequencies. With condition number tolerance chosen as \(\varepsilon_\mathrm{cond}=10^{48}\), the compact formula \eqref{eq:compact-formula} still provides a reasonable accuracy for a large range of frequencies to the right of the line. It is worth noting that it was found that the overestimation of the condition number is more severe for the polynomial system. The true condition number is about \(\kappa(Z(s,p))\approx10^9\) on the yellow line in \cref{fig:poly-2derr}.

The results are in general extremely good considering the limited amount of data used for the interpolation. Even in the case of polynomial dependence on the parameter, the algorithm proved itself to work well.

\begin{figure}
\centering
\includegraphics{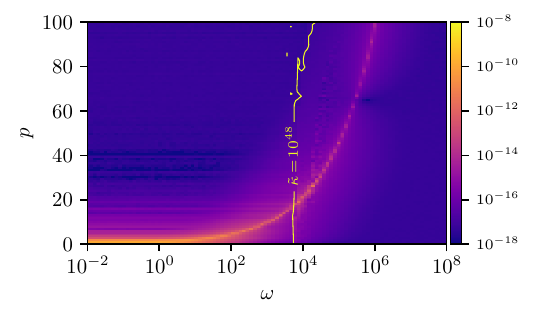}
\caption{Absolute error \eqref{eq:error} of the interpolated realizations over frequency and parameter for the polynomial example (\cref{ex:polynomial-num}).}
\label{fig:poly-2derr}
\end{figure}

\end{example}

\begin{example}[Parametric Penzl system]
\label{ex:penzl}
For the means of showcasing our method for larger systems, we consider the so-called Penzl system that has been extended to include an affine parameter dependency in \cite{ionita2014}. The model is given by the following realization:
\begin{align*}
A(p)&\!=\!\operatorname{diag}(A_1(p),A_2,A_3,A_4),\\
A_1(p)&\!=\!\begin{bmatrix*}[r]-1&p\\-p&-1\end{bmatrix*},\,A_2\!=\!\begin{bmatrix*}[r]-1&200\\-200&-1\end{bmatrix*},\,A_3\!=\!\begin{bmatrix*}[r]-1&400\\-400&-1\end{bmatrix*},\\
A_4&\!=\!-\operatorname{diag}(1,\,\dots,\,1000),\\
B^T&\!=\!\bigl[~\underbrace{10~\cdots~10}_{6}~\underbrace{1~\cdots~1}_{1000}~\bigr]=C
\end{align*}
The parametric transfer function $\hat{H}(s,p)$ is constructed based on only four uniformly distributed parameter values (the same as in \cref{ex:toy-num}) over the parameter interval $P=[0,100]$. As in the previous examples, we follow the steps of \cref{alg:1} with \(\varepsilon=10^{-7}\) and an alternating partitioning. The resulting Loewner matrices \(\mathbb{L},\mathbb{L}_s\in\mathbb{R}^{2014\times2014}\) are truncated to a rank of \(r=1009\) based on \eqref{eq:loewner-tol}.

Similarly to the other examples, the error plot in \cref{fig:penzl-2derr} shows the absolute error between the interpolated transfer function $\hat{H}(s,p)$ and the true values of the transfer function $H(s,p)$ according to \eqref{eq:error} with respect to different frequency and parameter values. The overestimation of the condition number is also present here, but less pronounced as in for \cref{ex:polynomial-num}. The true condition number is about \(\kappa(Z(s,p))\approx10^4\) on the yellow line in \cref{fig:penzl-2derr}.

\begin{figure}
\centering
\includegraphics{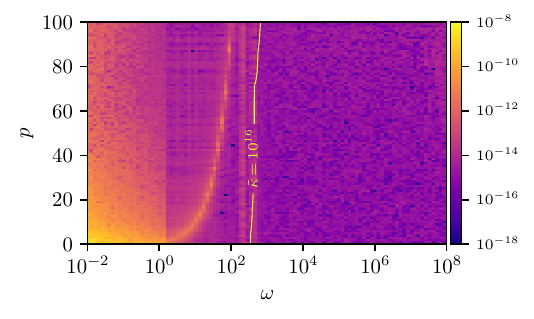}
\caption{Absolute error \eqref{eq:error} of the interpolated realizations over frequency and parameter for the parametric Penzl system (\cref{ex:penzl}).}
\label{fig:penzl-2derr}
\end{figure}

The results show that our method also works for large-scale systems where we are able to obtain a global parametric transfer function that also has a low error compared to the true transfer function for a large range of parameter values.
\end{example}

\section{Conclusion}
In this paper, we presented a new method to compute global parametric transfer function realizations based only on snapshot data. The method employs the univariate Loewner framework to interpolate the data. The proposed method is able to accurately generate a global parametric transfer function that is demonstrated to work for a wide range of frequencies and parameters. In the case of parametric polynomial dependence, rank bounds are provided for the Loewner pencil that interpolates the data and reveals the minimal order of the global parametric realization. The method is shown to work for different model classes and also for large-scale systems with only limited data. As future work, one could consider an adaptive sampling scheme in the parameter space to extract the most information from the data.

\section*{Acknowledgements}
The authors would like to thank Volker Mehrmann and Peter Benner for reading an early version of the manuscript and their valuable and insightful comments.

Art J. R. Pelling was supported by the German Research Foundation (DFG) in the project \emph{"Reduced order modelling of acoustical systems based on measurement data"} (\href{https://gepris.dfg.de/gepris/projekt/504367810?language=en}{504367810}).

\bibliographystyle{siamplain}
\footnotesize
\bibliography{references.bib}

\end{document}